\documentclass[a4paper,11pt]{article}
\usepackage[english]{babel}
\usepackage{amsfonts,amssymb,amsmath,amsthm,latexsym}

\newtheorem{rmk}{Remark}

\newtheorem{ex}{Example}

\newtheorem{theo}{Theorem}
\newtheorem{prop}{Proposition}

\newcommand{\C}{\mathbb{C}}

\newcommand{\fin}{\hfill$\square$\\}

\newcommand{\tx}{\textrm}

 \title{Plane branches with Newton nondegenerate polars}

\author{A. Hefez (UFF), M.E. Hernandes (UEM) and \\ M.F. Hern\'andez Iglesias (UEM)\thanks{The first two authors were partially supported by CNPq Grants and the third author by a fellowship from CAPES and ARAUCARIA foundations 
} }

\begin{document}

\maketitle
\begin{abstract}
We characterize the equisingularity classes of irreducible plane curve germs whose general members have a Newton nondegenerate general polar curve. In addition, we give explicit Zariski open sets of curves in such equisingularity classes whose general polars are Newton nondegenerate and describe their topology.
\end{abstract}

Authors' \, e-mail addresses: \, hefez@mat.uff.br, \ \ \ mehernandes@uem.br, \\
\centerline{ mfhiglesias.pma@uem.br}
\section{Introduction}

Let $\C\{x,y\}$ be the ring of convergent power series at the origin of $\C^2$. A germ of an analytic plane curve at the origin of $\C^2$ is a germ of set determined by an equation $f=0$, where $f$ belongs to the maximal ideal of the ring $\C\{x,y\}$ and has no multiple factors. Such a germ may be identified with the class $(f)$ of associated elements to the power series $f$, that is, the set of products of $f$ with the units in $\C\{x,y\}$. 

The power series $f$ will be called an equation of the germ of curve $(f)$. The multiplicity of $(f)$ is the multiplicity of its equation $f$ and $(f)$ will be said singular if its multiplicity is greater or equal than $2$. When $f$ is irreducible, we call $(f)$ a branch. 

We will say that two germs of analytic plane curves $(f)$ and $(f')$ are {\em analytically equivalent}, or equivalent for short, if there are a germ $\varphi$ of an analytic isomorphism at the origin of $\C^2$ and a unit $u$ of $\C\{x,y\}$ such that $f'=u f\circ \varphi$. We will say that $(f)$ and $(f')$ are {\em topologically equivalent}, or {\em equisingular}, if $\varphi$ is only a local homeomorphism. To verify that two germs of plane curves are equivalent is a hard task, while to verify that they are equisingular is rather easy, thanks to the classical work of Brauner and Zariski that we summarize below.  

Each branch of multiplicity $n$ is analytically representable by a Weierstrass polynomial $f$ in $\C\{x\}[y]$ of degree $n$ and there is a fractional power  series of the form $$
y=\varphi(x^\frac 1n)= x^{\frac mn}+\sum_{i\geq 1} a_i x^\frac{m+i}{n}, \quad \text{with} \quad m>n,$$
such that $f(x,\varphi(x^\frac 1n))=0$.

If $w$ is a primitive $n$-th root of unity, the conjugates $\varphi_j(x^{\frac 1n})=\varphi(w^{j}x^{\frac 1n})$, $1\leq j \leq n$, of $\varphi(x^{\frac 1n})$, will be called the Newton-Puiseux parametrizations of $(f)$. One has that $$f=\prod_{j=1}^{n}(y- \varphi_j(x^{\frac{1}{n}})).$$

The equisingularity class of a branch is encoded by the characteristic exponents of any of its Newton-Puiseux parametrization, which, in turn, determine and are determined by the set of intersection numbers of $f$ with all series in $\C\{x,y\}$ that are not multiple of $f$, called the {\em semigroup of values} of the branch (cf. \cite{Za} for all these notions). This semigroup will be represented by $\langle v_0,\ldots, v_g\rangle$ where the $v_i$'s are its minimal set of generators. The number $g$ is called the {\em genus} of the semigroup and is equal to the number of characteristic exponents.

More generally, the equisingularity class of a curve is encoded by the characteristic exponents of its branches and their mutual intersection numbers, where the intersection number of two curves $(f)$ and $(h)$ is defined as 
$$
{\rm I}(f,h)=\dim_\C \frac{\C\{x,y\}}{\langle f,h \rangle},
$$
where $\langle f,h \rangle$ denotes the ideal generated by $f$ and $h$ in the ring $\C\{x,y\}$ and the dimension we consider is as a $\C$-vector space.

If $f=\sum a_{i,j}x^{i}y^j\in \C\{x,y\}$, we will denote by $N(f)$ the Newton polygon of $f$ or of $(f)$, since the Newton polygon is invariant by multiplication by units. If $\ell$ is a side of $N(f)$, we define 
$$f_\ell=\sum_{(i,j)\in \ell} a_{i,j}x^{i}y^j.$$

We will say that a reduced $f$ is {\em Newton nondegenerate}, or simply nondegenerate, with respect to the coordinate system  $(x,y)$, if for every side $\ell$ of $N(f)$, the polynomial $f_\ell$ has no critical point outside the lines $x=0$ and $y=0$. This is equivalent to say that the polynomial
$F_\ell(z)=\frac{f_\ell(1,z)}{z^{j_0}}$, where $j_0=\min \{ j ; \ (i,j)\in \ell \}$, has no multiple roots. We call $F_\ell(z)$ the {\em associated polynomial} to the side $\ell$.

Notice that the notion of nondegenerate is also invariant by multiplication by units. So, to say that $f$ is nondegenerate is the same as to say that the germ of curve $(f)$ is nondegenerate.\medskip

The importance of the notion of nondegenerate curves comes from the fact that the topology of such curves is completely determined by their Newton polygons as follows from the result below (cf. \cite[Propostion 4.7]{O}).\medskip

\noindent {\bf Oka's Lemma}
 {\em Let $(f)$ be a reduced plane curve germ and $(x,y)$ be a local coordinate system such that the tangent cone to $(f)$ does not contain the line $(x)$. If $(f)$ is nondegenerate, then for each side $\ell_k$ of $N(f)$ with extremal points $(i_k,j_k)$ and $(i_k',j_k')$, with $j_k' > j_k$ and $i_k > i_k'$, there correspond $d_k=\tx{GCD}(n_k,m_k)$ branches $\xi_{k,l}$, $1\leq l \leq d_k$, where  $n_k=j_k'-j_k$ and $m_k=i_k-i_k'$, with semigroup generated by $\frac{n_k}{d_k}$ and $\frac{m_k}{d_k}$ (allowing any one of them to be $1$). If $\xi$ and $\xi'$ are two of these branches with semigroups generated by $\alpha_0,\alpha_1$ and $\beta_0,\beta_1$, respectively, then their intersection number is given by
 	$ {\rm I}(\xi,\xi')=\min\{\alpha_0\beta_1,\beta_0\alpha_1\}.$}\medskip

Oka's Lemma shows that a nondegenerate curve may only have nonsingular branches or branches of genus one.

As a consequence, it follows that two reduced nondegenerate power series with same Newton polygons determine equisingular curves.

Let us define now our main object of study. The {\em general polar curve} of $f$ is the curve determined by an equation $P(f):af_x+bf_y=0$, where $(a:b)$ is a general point in $\mathbb{P}^1_{\mathbb{C}}$.

It is known that for equivalent $f$ and $f'$, their general polar curves are  equisingular, while two equisingular $(f)$ and $(f')$ have not necessarily equisingular general polars (cf. \cite{T} and \cite{P}).

We will now focus on the study of equisingularity classes of plane branches whose general members have Newton nondegenerate polar curves. One of the known general results about polars is  \cite[Theorem 3.1]{merle}, which in particular implies that the general polar curve of a plane branch of genus greater than two always has, at least, one branch of genus greater than 1, and consequently is Newton degenerate. Hence, to study non-degeneracy conditions for $P(f)$ it is sufficient to consider $(f)$ of genus one or two. In Section 2 we will focus on the case of genus $1$, where the content of \cite{C1} is summarized and condensed through the use of the concept of Newton nondegenerate curves. These results will be used in Section 3 to describe the topology of the general plane branches of genus $2$ with nondegenerate general polars.

\section{Nondegenerate polars of branches of genus $1$}

It is known (cf. \cite{Za} ) that an irreducible plane curve has a semigroup $\langle p,q \rangle$, with ${\rm GCD}(p,q)=1$, if and only if it has an analytic representative with an equation of the form 
$$f_1=y^p-x^q+\sum_{ip+jq > pq}a_{i,j}x^iy^j, \ \ a_{i,j}\in \C.$$

In what follows, we denote by $K(p,q)$ the set of all these equations. 
In this way, the polar curve of an element $f_1\in K(p,q)$ is given
by an equation 

$$P(f_1)=-qax^{q-1}+pby^{p-1}+ \sum_{ip+jq > pq}
(iaa_{i,j}y+ 
jba_{i,j}x)x^{i-1}y^{j-1}.$$

Let $f_1$ be a general member of $K(p,q)$. For each fixed $j$, it is easy to verify that the lowest exponent of $x$ in a term of $(f_1)_y$ having as factor the power $y^j$ is $\alpha(j) = q-\left [\frac{(j+1)q}{p}\right ]$, where $[k]$ represents the integer part of $k$. On the other hand, the lowest exponent of $x$ in a term of $(f_1)_x$ having as factor the power $y^j$ is $\gamma(j) = q-\left(\left [\frac{jq}{p}\right ] +1 \right)$. Since $\alpha(j)\leq \gamma(j)$, it follows that the powers of the monomials that determine $N(P(f_1))$ are contained in the set
\begin{equation}
\label{pontos} A= \left \{(\alpha(j), j); \ 0 \leq
j \leq p-1 \right\} \end{equation} and each element
$(\alpha(j),j)\in A$ is associated to the following term of $P(f_1)$:\medskip

\noindent $ t_j=(j+1)ba_{\alpha(j),j+1}x^{\alpha(j)}y^j, \qquad \tx{if}\,\,
    \left[\frac{(j+1)q}{p}\right]\neq \left[\frac{jq}{p}\right]+1, \quad \text{or}$\medskip
    
\noindent $
t_j=\big( b(j+1)a_{\alpha(j),j+1}+a(q-[(jq/p])a_{\alpha(j)+1,j} \big)x^{\alpha(j)}y^j
    ,\,\,\,  \tx{if}\,\,\left[\frac{(j+1)q}{p}\right]= \left[\frac{j q}{p}\right]+1.
        $\medskip

Let $\frac{q}{p}=[h_0,\ldots,h_s]$ be the continued fraction decomposition of $\frac{q}{p}$ and consider, for $i\leq s$, $\frac{q_i}{p_i}=[h_0,\ldots, h_i]$ (here necessarily  $\tx{GCD}(p_i,q_i)=1$),

In \cite{C1} it is shown that the Newton polygon $N(P(f_1))$, determined by the set $A$ in 
(\ref{pontos}), has sides $\ell_k$, where $0\leq k \leq \frac{s-1}{2}$ (with $s$ as above). We then have that $\ell_k$ contains the points
 \begin{equation} \label{pointsgenusone}
\begin{array}{ll} 
 (q-q_{2k},p_{2k}-1) +  l(-q_{2k+1},p_{2k+1}); \,\, 0 \leq l \leq h_{2k+2},  & {\rm when} \ 0 \leq k < \frac{s-1}{2},  \\ \\
 (q-q_{s-1},p_{s-1}-1) \ \text{and} \ (0,p-1),  & {\rm when} \ k=\frac{s-1}{2},
\end{array}
\end{equation}
where the last possibility only occurs when $s$ is odd.

Using (\ref{pointsgenusone}) and considering the terms $t_j$, we have that the associated polynomial $F_k(z)$ to the side $\ell_k$ of the polygon $N(P(f_1))$ is the following:
 $$
F_k(z)=\sum_{l=0}^{h_{2k+2}}\frac{t_{l
        p_{2k+1}+p_{2k}-1}(1,z)}{z^{p_{2k}-1}}, \quad  \ 0\leq k < \frac{s-1}{2} $$
and when $s$ is odd
$$
F_{\frac{s-1}{2}}(z)= \frac{t_{p_{s-1}-1}(1,z)+t_{p-1}(1,z)}{z^{p_{s-1}-1}}.
$$

For a polynomial $F(z)$ we denote its
discriminant and its zero set by $\Delta(F(z))$ and $Z(F(z))$, respectively.\smallskip

Let us consider the Zariski closed set $Z(p,q)$ of elements in $K(p,q)$ whose coefficients are in the set  \begin{equation} \label{zeros}
 Z \left( \prod_\lambda t_\lambda(1,1) \prod_\mu\Delta(F_\mu(z)) \right),
 \end{equation}
where $\lambda$ varies in the set $\{p_{2k}-1+lp_{2k+1}; \ 0\leq k < \frac{s-1}{2}, \ 0\leq l \leq h_{2k+2}\}$, union with the set $\{p_{s-1}-1, p-1\}$ if $s$ is odd; and $0\leq \mu \leq \frac{s-1}{2}$.

With the above notation we have the following proposition.

 \begin{prop}\label{genusone}
 If $f_1\in K(p,q)\setminus Z(p,q)$, then the polar curve $P(f_1)$ is Newton nondegenerate.
 \end{prop}

\noindent {\bf Proof} \ The result is clear since the condition $f_1\not \in Z(p,q)$ implies that $N(P(f_1))$ supports exactly the points in (\ref{pointsgenusone}) corresponding to monomials of the polar curve and that the associated polynomials to the side of this polygon have no multiple roots.
\fin

Combining Proposition and Oka's Lemma, we get the following result.
\begin{theo}
\label{topgenus1} The polar $P(f_1)$ of an element $f_1\in
K(p,q)\setminus Z(p,q)$ has branches $\xi_{k,\lambda}$  for \ $0 \leq k < \frac{s-1}{2}$ and $1 \leq  \lambda \leq h_{2k}$ with semigroup $\langle p_{2k+1},q_{2k+1} \rangle$. If $s$ is odd, one has additionally $d_s=\tx{GCD}(p_s-p_{s-1}, q_s-q_{s-1})$ branches with semigroup $\langle \frac{p_s-p_{s-1}}{d_s}, \frac{q_s-q_{s-1}}{d_s} \rangle$. Moreover, if $\xi$ and $\xi'$ are branches of $P(f_1)$ with semigroups $\langle v_0,v_1\rangle$ and $\langle w_0,w_1\rangle$, respectively,  then $I(\xi,\xi')=min\{v_0w_1,w_0v_1\}$ and  if $\xi$ is smooth then $I(\xi,\xi')=w_1$.
\end{theo}

\begin{ex}\rm 
    Let $f_1=y^7-x^{19}+\sum_{7i+19j > 133}a_{i,j}x^iy^j$. We have that $\frac{19}{7}=[2,1,2,2]=[h_0,h_1,h_2,h_3]$.  Since $h_2=2$, corresponding to the side $\ell_0$, the polar  $P(f_1)$ has $2$ branches $\xi_1$ and $\xi_2$, which are smooth because $ p_1=1$. Since $s=3$ is odd, $p_{s-1}=3$, $q_{s-1}=8$, $p_s=7$ and $q_s=19$, we have that $1=\tx{GCD}(4,11)=\tx{GCD}(p-p_{s-1},q-q_{s-1})$, then $P(f_1)$ has an extra branch $\rho$ with semigroup $\langle 4, 11 \rangle$, corresponding to the side $\ell_1$.

The side $l_0$ is determined by the terms $t_{0}=ba_{17,1}x^{17}$,
    $t_{1}=2ba_{14,2}x^{14}y$ and  $t_{2}=3ba_{11,3}x^{11}y^2$.  In
    this way, we have $F_0(z)=3ba_{11,3}z^2+2ba_{14,2}z+ba_{17,1}$ and
    $\Delta(F_0(z))=12b^3a_{11,3}(3a_{11,3}a_{17,1}-a_{14,2}^2)$. The side $l_1$
    is determined by $t_{2}=3ba_{11,3}x^{11}y^{2}$ and $t_{6}=7by^6$, so
    $F_1(z)=b(7z^4+3a_{11,3})$ and $\Delta(F_1(z))=4(84b^2a_{11,3})^3$.

    Therefore, the set (\ref{zeros}) given before Proposition \ref{genusone}
    is 
   $$Z(a_{11,3}a_{14,2}a_{17,1}(3a_{11,3}a_{17,1}-a_{14,2}^2)).$$ 

So, $Z(7,19)$ is the subset of $K(7,19)$ whose coefficients belong to the above set, hence 
    $K(7,19)\setminus Z(7,19)$ is an open set where the curves $(f_1)$ have  Newton
    nondegenerate general polar curves. Finally one has that $I(\xi_1,\xi_2)=3$ and $I(\xi_i,\rho)=11$ for $i=1,2$.
 \end{ex}

Notice that, for all $f_1\in K(p,q)\setminus Z(p,q)$, the branches of the general polar curve have at most genus one. For $f_1 \in Z(p,q)$ this may fail as we show in the next example.

 \begin{ex} \rm If $f_1\in K(5,12)$, then the set (\ref{zeros}) described in Proposition \ref{genusone} is
$$
Z(a_{5,3}a_{10,1}(9a_{5,3}^2-20a_{10,1})).
$$

Let $f_1=y^5-x^{12}+x^5 y^3+x^8y^2+\frac{9}{20}x^{10}y \ \in Z(5,12)$. Its polar curve is given by the polynomial
$$
P(f_1) = 5by^4 +5ax^4y^3+(8ax^7+3bx^5)y^2+(\frac{9}{2}ax^9+ 2x^8)y-12ax^{11}+\frac{9}{20}bx^{10}.
$$ By a computation using $\text{Maple}^{\text{TM}}$, we find that $P(f_1)$
 has a single branch with a parametrization given by
 $$ x=t^4 ,\,\,\,\,y=d_1t^{10}+d_2t^{11}+\cdots, \ \ \text{with} \ d_1d_2\neq 0.$$ Therefore, $P(f_1)$ is a branch with characteristic exponents $4,10,11$, hence it has genus two.
In \cite{polar} we give the complete description of the polars of curves in $K(5,12)$.
\end{ex}

The above example shows that the question of which equisingularity classes have all curves with nondegenerate general polars is meaningless. This is why we asked the question only for general members of equisingularity classes.  

\section{Nondegenerate polars of branches of genus $2$}

Any irreducible curve $(f)$ of genus 2 has a semigroup of values of the form $$\langle
v_0,v_1,v_2\rangle=\langle
e_1p,e_1q,e_1pq+d\rangle,$$ 
where $e_1={\rm GCD}(v_0,v_1)$ and ${\rm GCD}(e_1,d)=1$.

Using the construction of \cite{AB}, one knows that any $(f)$ with semigroup as above is equivalent to one of the form $f=f_1^{e_1}+f_2$, where
$f_1=y^{p}-x^q+\sum_{ip+jq > pq} a_{i,j}x^iy^j\in K(p,q)$ and 
 $$f_2=\sum_{iv_0+jv_1 > 2v_2}a_{i,j,2}x^iy^jf_1^{e_1-2}+\cdots+
\sum_{iv_0+jv_1 > (e_1-1)v_2}a_{i,j,e_1-1}x^iy^jf_1+$$
$$+b_{i_0,j_0}x^{i_0} y^{j_0}+\sum_{iv_0+jv_1> e_1v_2}b_{i,j}x^iy^j,$$
with $i_0v_0+j_0v_1=e_1v_2, $  \ $i_0\geq 0$, $0\leq j_0<\frac{v_0}{e_1}$; \ $a_{i,j},a_{i,j,k},b_{i,j}\in
\C$  \,\, and \,\, $b_{i_0,j_0}\neq 0.$ We denote by $K(e_1p,e_1q,e_1pq+d)$ the set of the   power series described above.

Notice that
$P(f)=e_1f_1^{e_1-1}P(f_1)+ P(f_2)$.

With an elementary computation involving the inequalities relating the indices in the summations in the expression of $f_2$, one may conclude that $N(f)=N(f_1^{e_1})$ and since $f_1$ is irreducible, we have that $N(f_1^{e_1})=N((y^p-x^q)^{e_1})$, therefore $N(f)=N((y^p-x^q)^{e_1})$, which is precisely the segment $l$ contained in the line
$pX+qY=e_1pq$ supporting only the points $(iq,(e_1-i)p)$ for $i=0,\ldots ,e_1$. 

We denote by $l_{p,q}$ the segment that contains
the points $(iq,(e_1-i)p-1)$ for $i=0,\ldots ,e_1-1$. 

Now, since the segments $l$ and
 $l_{p,q}$ are parallel with slope of absolute value less than $1$ and the second coordinate of their points are shifted by one unit, it follows that the segment $l_{p,q}$ belongs to the polygon 
$N(P(f))$.

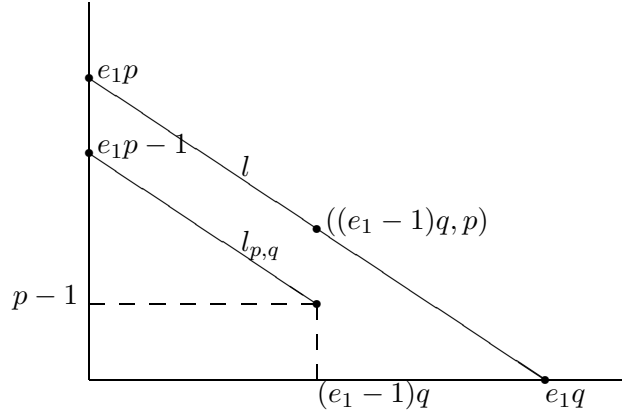
\begin{center}
    \setlength{\unitlength}{1cm}
    \begin{picture}(7,5)
    \put(0,0){\line(1,0){7}}\put(0,0){\line(0,1){5}}\put(2,1.7){$l_{p,q}$}
    \put(0,3){\line(3,-2){3}}
    \put(0,3){\circle*{0.1}$e_1p-1$}
    \put(0,4){\line(3,-2){6}}\put(0,4){\circle*{0.1}$e_1p$}
    \put(3,1){\circle*{0.1}}\put(3,2){\circle*{0.1}$((e_1-1)q,p)$}
    \put(0,0){\dashbox{0.2}(3,1){}}\put(-1,1){$p-1$}\put(3,-0.3){$(e_1-1)q$}
    \put(6,0){\circle*{0.1}}
    \put(6,-0.3){$e_1q$}
    \put(2,2.7){$l$}
    \end{picture}
 \end{center}\medskip

From the expression of $f$ we have that the terms of $P(f)$
corresponding to the points on $l_{p,q} $ are obtained from $e_1f_1^{e_1-1}P(f_1)$ and the associated polynomial to
$l_{p,q}$ is (modulo a nonzero constant) $F_{p,q}(z)=(z^p-1)^{e_1-1}$.
Therefore, for $e_1 > 2$ we have that $P(f)$ is a degenerate curve, and since we are interested in branches with nondegenerate general polars, it will be sufficient to consider $e_1=2$, that is, the case in which $f$ has a semigroup of the form $\langle 2p,2q,2pq+d \rangle$, with odd $d$, which are coincidentally the equisingularity classes studied in \cite{LP} in view of their analytic classification.

So, we are reduced to study polars of curves in the sets of the form $K(2p,2q,2pq+d)$, that is,  polar curves of series of the form $f=(f_1)^2+f_2$, where
$$f_1\in K(p,q), \, \, \, \, f_2= b_{i_0,j_0}x^{i_0} y^{j_0}+\sum_{ip+jq> 2pq+d}b_{i,j}x^iy^j,$$
where \ $i_0p+j_0q=2pq+d$, \ $i_0\geq 0,\  0\leq j_0 < p.$

As we have observed above, in this situation, one has that 
\begin{equation}\label{polargen2}
P(f)=2f_1P(f_1)+P(f_2)
\end{equation} and that the Newton polygon of $f$ is supported by the line $L:pX+qY=2pq$ and it
contains the points $(0,2p),\ (q,p)$ and $(2q,0)$. On the other hand, the segment $l_{p,q}$ (with $e_1=2$) with extremal points $(0,2p-1)$ and
$(q,p-1)$ is contained in one side of the polygon $N(P(f))$.

Notice that the Newton polygon of $f_1(P(f_1))$ is the Minkowski sum of $N(f_1)$ and $N(P(f_1))$; that is, if $A$ is as in (\ref{pontos}), then 
$N(f_1(P(f_1)))$ is determined by
$$\{(0,p-1)+A, (q,0)+A\}.$$ 

Since the absolute value of the slopes of the sides of $N(P(f_1))$ are
$\frac{p_{2j+1}}{q_{2j+1}}$ which are smaller than $\frac{p}{q}$ (the absolute value of the 
slope of $N(f_1)$), the polygon $N(f_1(P(f_1)))$ is
determined by the set $\{(0,2p-1), (q,0)+A\}$.

Therefore, $N(f_1(P(f_1)))$ is
given by $l_{p,q}$ and $(0,q)+N(P(f_1))$, which implies that $N(f_1P(f_1))$ has $[\frac{s+1}{2}]+1$ sides $D_k$ where
     \begin{equation}\label{puntos}
        D_k=(q,0)+l_k ,\,\,\,  0\leq k \leq  \frac{s-1}{2}\ \ \ \text{and} \ \ D_{[\frac{s-1}{2}]+1}=l_{p,q},
     \end{equation}
where the $l_k$ are given as in Section 2 and the
terms of $f_1P(f_1)$ associated to $(i,j)\in N(f_1P(f_1))$ is
$$
g_{j}=
   \begin{cases}
    -x^qt_{j}\,\,\,&  \text{if}\,\,\, j\neq 2p-1\\ py^{2p-1}  & \text{if} \,\, j= 2p-1.
    \end{cases}
$$

By a similar analysis to the case of genus one,
considering the polar $P(f_2)=a(f_2)_x+b(f_2)_y$, one has 
that for each $j$ the smallest value of $i$ for which there exists a term in $b(f_2)_y$ containing $y^j$ as a factor is $2q-[\frac{q(j+1)-d}{p}]$ and in $a(f_2)_x$ is 
$2q-[\frac{qj+p-d}{p}]$. Denoting by $\beta(j)$ the minimum of these values $2q-[\frac{q(j+1)-d}{p}]$, the term of 
$P(f_2)$ that contains the monomial $x^{\beta(j)}y^j$ 	is\medskip

\noindent $
     h_{j}=     (j+1)bb_{\beta(j),j+1}x^{\beta(j)}y^j,\ \ \tx{if}\,\,
        \beta(j)\neq 2q-[\frac{jq+p-d}{p}];  \ \
\tx{or}$\medskip

\noindent $h_j= \big( b(j+1)b_{\beta(j),j+1}+a(2q-[(jq+p-d)/p]+1)b_{\beta(j)+1,j}\big)x^{\beta(j)}y^j,$
\medskip

\noindent if \ \ $ \beta(j) = 2q-[\frac{jq+p-d}{p}]$.\medskip
          
Since $\beta(j)\geq \alpha(j)+q$, we conclude that the polygons
$N(P(f))$ and $N(f_1P(f_1))$ are equal and have the same
points on the respective sides. So, from Formula (\ref{polargen2}), it follows that the term containing $x^{\alpha(j)+q}y^j$ in $P(f)$ is
$$
 m_{j}=
     \begin{cases}
    2g_{j}+h_{j},\,\,& \tx{if}\,\,
    \alpha(j) +q=\beta(j)\\
    2g_{j} ,&  \tx{if}\,\,\alpha(j) +q < \beta(j).
    \end{cases}
 $$
 
From the above expression of the terms $m_{j}$ and from (\ref{puntos}), the associated polynomial to the side $D_k$ of
$N(P(f))$ is 
$$ \begin{array}{l}
F_k(z)=\sum_{l=0}^{h_{2k+2}}\frac{m_{p_{2k}-1+lp_{2k+1}}(1,z)}{z^{p_{2k}-1}}, \ \ \tx{if} \ \ 0\leq k < \frac{s-1}{2},\\
F_{\frac{s-1}{2}}(z)= \frac{m_{p_{s-1}-1}(1,z)+m_{p-1}(1,z)}{z^{p_{s-1}-1}} \ \ \quad \tx{and} \\
F_{[\frac{s-1}{2}]+1}(z)=\frac{m_{2p-1}(1,z)+m_{p-1}(1,z)}{z^{p-1}}.
\end{array}
$$

Let us consider the Zariski closed set $Z(2p,2q,2pq+d)$ of elements in $K(2p,2q,2pq+d)$ whose coefficients are in the set  \begin{equation}
Z \left( \prod_\lambda m_\lambda(1,1) \prod_\mu\Delta(F_\mu(z)) \right),
\label{zeros2}
\end{equation}
where $\lambda$ varies in the set $\{p_{2k}-1+lp_{2k+1}; \ 0\leq k < \frac{s-1}{2}, \ 0\leq l \leq h_{2k+2}\}\cup\{2p-1\}$, union with the set $\{p_{s-1}-1, p-1\}$ if $s$ is odd; and $0\leq \mu \leq [\frac{s+1}{2}]+1$.

By the above analysis and Oka's Lemma we obtain the following
result.
\begin{theo}\label{genus2}
    If $f \in K(2p,2q,2pq+d)\setminus Z(2p,2q,2pq+d)$, then the polar curve $P(f)$
    is Newton nondegenerate. In this case, $P(f)$ has a branch with semigroup
    $\langle p,q\rangle$ and $\left [
    \frac{s+1}{2}\right ]$ branches as described in Theorem \ref{topgenus1}, such that,
    if $(\xi)$ and $(\xi')$ are branches of $P(f)$ with semigroup
    $\langle v_0,v_1\rangle$ and $\langle w_0,w_1\rangle$, then
    $I(\xi,\xi')=min\{v_0w_1,w_0v_1\}$, if $(\xi)$ is smooth, then $I(\xi,\xi')=w_1$.
\end{theo}

\begin{rmk} In \cite{polar} we showed that, for any $f\in K(4,2k,4k+d)$, with $k>1$ odd, the polar
	curve $P(f)$ has a smooth branch $\xi$ and a branch $\rho$ with
	semigroup $\langle 2,k\rangle$ such that $I(\xi,\rho)=k$ and $Z(4,2k,4k+d)$ is empty.
	\end{rmk}

\begin{ex} \rm If $f\in K(10,24,121)$ is given by
    $$f=(y^5-x^{12}+\sum_{5i+12j > 60}a_{i,j}x^iy^j)^2+b_{17,3}x^{17}y^{3}
    +\sum_{5i+12j > 121}b_{i,j}x^iy^j,$$ with $b_{17,3}\neq 0$, then the Newton polygon
    $N(P(f))$ is given by the following diagram where we wrote the term $m_j$ at the point with coordinates $(\beta(j),j)$\bigskip
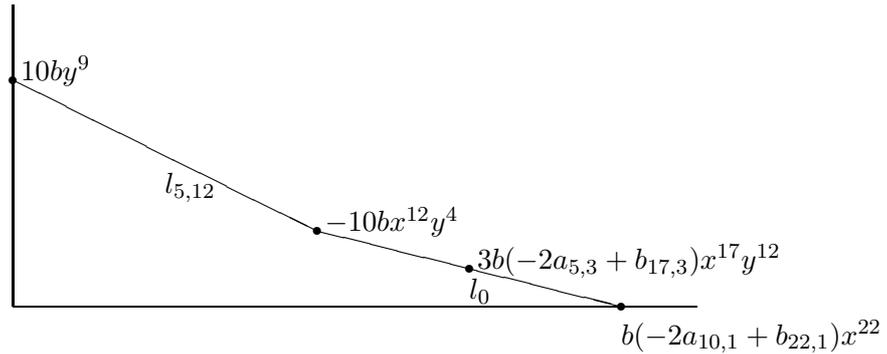
\begin{center}
    \setlength{\unitlength}{1cm}
    \begin{picture}(9,4)
    \put(0,0){\line(1,0){9}}\put(0,0){\line(0,1){4}}
    \put(0,3){\line(2,-1){4}}
    \put(4,1){\line(4,-1){4}}
    \put(0,3){\circle*{0.1}$10by^9$}
    \put(4,1){\circle*{0.1}$-10bx^{12}y^4$}
    \put(8,0){\circle*{0.1}}
    \put(8,-0.5){$b(-2a_{10,1}+b_{22,1})x^{22}$}
    \put(6,0.5){\circle*{0.1}$3b(-2a_{5,3}+b_{17,3})x^{17}y^{12}$}

    \put(2,1.5){$l_{5,12}$}
    \put(6,0.1){$l_0$}
    \end{picture}
\end{center}\vspace{.8cm}
and the associated polynomial to the side $l_0$ is
$F_0(z)= b(-10z^4 + 3(b_{17,3}-2a_{5,3})z^2 + b_{22,1}-2a_{10,1})$.

In this case, we have that the set in (\ref{zeros2}) is
$$Z((b_{17,3}-2a_{5,3})(b_{22,1}-2a_{10,1})(9b_{17,3}^2-36a_{5,3}b_{17,3}+36a_{5,3}^2+40b_{22,1}-80a_{10,1})).$$

\end{ex}

\begin{rmk}\label{dcondition}
If $f=f_1^2+f_2\in K(2p,2q,2pq+d)$ and $d\geq q$, then the set
$Z(2p,2q,2pq+d)$ depends only on $f_1$. In fact, since 
$f_2=\sum_{ip+jq\geq 2pq+d}b_{ij}x^iy^j$, if $d\geq q$, we have
$ip+jq\geq 2pq+d\geq 2pq+q > 2pq+p$ and in this way all vertices of
$N(P(f_2))$ are above $N(f_1^2)$. Hence the monomial
$m_{\alpha(j),j}$ does not depend on the coefficients of
$f_2$.

\end{rmk}

\begin{rmk}
It is known that the topological class of the general polar curve contains analytical data of the curve, but does not allow to recover the topological class of the curve. For instance, if we consider general elements $f_k\in K(2p,2q,2pq+d_k)$, then by the previous remark $P(f_k)$ belongs to the same topological class for all $d_k\geq q$, but the topological class of $f_k$ varies according to $d_k$. In contrast, in \cite{C2} it was shown that the topological type of a curve may be recovered from the cluster of basis points of the general polar curve. 
\end{rmk}

We summarize the above discussion in the following statement.

\begin{theo}
 The general polar curve $P(f)$ of the general member of an equisingularity class of branches is Newton nondegenerate if and only if the equisingularity class corresponds to one of the following semigroups: $\langle p,q \rangle $ or $\langle 2p,2q,2pq+d\rangle $, where GCD$(p,q)=1$ and $d$ is odd. The topology of the general polar in each case is described in Theorem \ref{genusone} and in Theorem \ref{genus2}.
 \end{theo}

\end{document}